\documentclass{article}

\usepackage{amssymb}
\usepackage{amsmath}
\usepackage{mathrsfs}
\usepackage[all]{xy}

\newtheorem{thm}{Theorem}[section]
\newtheorem{pro}[thm]{Proposition}
\newtheorem{lem}[thm]{Lemma}
\newtheorem{cor}[thm]{Corollary}
\newtheorem{prb}[thm]{Problem}
\newtheorem{rem}[thm]{Remark}
\newtheorem{exa}[thm]{Example}

\def\F{\mathfrak{F}}

\def\dd{\mathbf{d}}
\def\w{\mathbf{w}}
\def\es{\varnothing}

\def\Ra{\Rightarrow}

\def\st{^\ast}

\def\ol#1{\overline{#1}}
\def\wh#1{\widehat{#1}}
\def\Fr{Fra\"{\i}ss\'e\ }
\def\C{\mathscr{C}}

\def\T{\mathcal{T}}
\def\cd{{\rm ($\dag$)}}

\def\beginproof{\noindent\textit{Proof.} }
\def\endproof{\hfill $\Box$ \bigskip}

\renewcommand\leq{\leqslant}
\renewcommand\geq{\geqslant}

\DeclareMathOperator\End{End}  \DeclareMathOperator\Sym{Sym} \DeclareMathOperator\krn{ker}
   
\DeclareMathOperator\Fl{Flim} 

\title%[The Bergman Property for Endomorphism Monoids]%
{The Bergman property for endomorphism monoids\\ of some Fra\"{\i}ss\'e limits}

\author{Igor Dolinka}

\date{}

\begin{document}

\maketitle

\begin{abstract}
Based on an idea of Y.~P\'eresse and some results of Maltcev, Mitchell and Ru\v skuc, we present sufficient conditions
under which the endomorphism monoid of a countably infinite ultrahomogeneous first-order structure has the Bergman
property. This property has played a prominent role both in the theory of infinite per\-mu\-ta\-tion groups and, more
recently, in semigroup theory. As a byproduct of our considerations, we establish a criterion for a countably infinite
ultrahomogeneous structure to be homomorphism-homogeneous.\medskip

\footnote{Supported by Grant No.174019 of the Ministry of Science and Technological Development of the Republic of
Serbia.}
\end{abstract}

\section{Introduction}

\subsection{The Bergman property and \Fr limits}

Let $S$ be a semigroup. For $\es\neq A\subseteq S$ and $n\geq 1$ we denote
$$A^n=\{a_1\cdots a_n:\ a_1,\dots,a_n\in A\}.$$
The least subsemigroup of $S$ containing $A$ is said to be \emph{generated} by $A$; usually it is denoted by $\langle A
\rangle$. Clearly,
$$\langle A \rangle=\bigcup_{n\geq 1}A^n.$$
In particular, if $S=\langle A\rangle$ then $A$ is a \emph{generating set} of $S$, thus by definition $S=\bigcup_{n\geq
1}A^n$. However, it might turn out that only a finite portion of the latter infinitary union suffices to obtain the
whole $S$, that is, $S=\bigcup_{n=1}^m A^n$ holds for some $m\geq 1$. In such a case we say that $S$ is \emph{semigroup
Cayley bounded} with respect to $A$. For a group $G$ and its (group) generating set $\Gamma$ we have that $G$ is
generated as a semigroup by $\Gamma\cup\Gamma^{-1}$ (where $\Gamma^{-1}=\{g^{-1}:\ g\in\Gamma\}$); thus we say that $G$
is \emph{group Cayley bounded} with respect to $\Gamma$ if it is semigroup Cayley bounded with respect to
$A=\Gamma\cup\Gamma^{-1}$ (that is, the Cayley graph of $G$ with respect to $\Gamma$ is of finite diameter). A
well-known result of George Bergman \cite{B} asserts that for any (infinite) set $X$, the symmetric group $\Sym(X)$ is
group Cayley bounded with respect to \emph{every} generating set. Hence, the term `the Bergman property' quickly
established itself \cite{DG,Kh} to describe the property of groups of being group Cayley bounded with respect to every
generating set. To distinguish between groups and semigroups, we refer to this remarkable property as the \emph{group}
Bergman property; the analogous property for semigroups---the main subject of investigation in a recent contribution by
Maltcev, Mitchell and Ru\v skuc \cite{MMR}---is called the \emph{semigroup} Bergman property. For a group $G$, the
semigroup Bergman property obviously implies the group Bergman property. It is still unknown, however, whether the
converse is true.

The other principal theme of this paper are the fascinating objects from model theory called \emph{\Fr limits}. Namely,
if $\C$ is a countable set of finitely generated first-order structures of a fixed countable signature that is closed
for taking (isomorphic copies of) substructures, has the \emph{joint embedding property} (JEP) and the
\emph{amalgamation property} (AP), then a celebrated result of Roland \Fr \cite{F1,H} guarantees the existence and
uniqueness of a countable structure $F$ such that:
\begin{itemize}
\item[(i)] the set of all finitely generated substructures of $F$, called the \emph{age} of $F$, coincides (up to
isomorphism) with $\C$, and
\item[(ii)] $F$ is \emph{ultrahomogeneous}, which means that for any isomorphism $\alpha:A\to A'$ between finitely
generated substructures $A,A'$ of $F$ there is an automorphism $\wh\alpha$ of $F$ that extends $\alpha$, i.e.\
$\wh\alpha|_A=\alpha$.
\end{itemize}
Following \cite{H}, such a structure $F$ is called the \emph{\Fr limit of} $\C$ and denoted by $\Fl(\C)$. Moreover, any
countably infinite ultrahomogeneous structure arises in this way: it is the \Fr limit of the class of all of its
finitely generated substructures. A class of finitely generated structures that satisfies the premises of the \Fr
theorem is called a \emph{\Fr class}. It is not difficult to show that if $\C$ is a \Fr class, and if $\ol\C$ denotes
the class of all countable structures whose ages consist of structures isomorphic to members of $\C$, then in fact any
member of $\ol\C$ embeds into $\Fl(\C)$. Historically, the first \Fr limits discovered were the rational Urysohn space
$\mathbb{U}_\mathbb{Q}$ \cite{Ur} as the limit of all finite metric spaces with rational distances, and $\mathbb{Q}$,
the limit of all finite chains \cite{F1}. Other classical examples of \Fr classes and their limits include:
\begin{itemize}
\item finite simple graphs and the \emph{random graph} $R$ \cite{C1,C2},
\item finite posets and the \emph{generic poset} $\mathbb{P}$ \cite{Sch},
\item finite semilattices and the countable universal ultrahomogeneous semilattice $\Omega$ \cite{DKT},
\item finite distributive lattices and the countable universal ultrahomogeneous distributive lattice $\mathbb{D}$
\cite{DrM},
\item finite Boolean algebras and the countable atomless Boolean algebra $\mathbb{A}$.
\end{itemize}

The paper \cite{MMR} provides an abundance of examples of well-known semigroups both with and without the Bergman
property. A significant part of those examples are semigroups of various mappings, or even morphisms of some structure.
However, it is one particular result contained in \cite[Theorem 4.2]{MMR} that will be of a special interest here: this
is the assertion that $\End(R)$, the endomorphism monoid of the random graph $R$, has the Bergman property. This claim
remained unproved in \cite{MMR}; the theorem itself was formulated as a consequence of Lemma 2.4 of that paper (see
Lemma \ref{sd} below) and two earlier publications \cite{AMS,MPQ}, which indeed account for all the assertions
contained in the theorem except for the one about $\End(R)$. Later I learned \cite{JDM} that the Bergman property for
$\End(R)$ is a consequence of a result in the recent doctoral thesis of Y.~P\'eresse \cite{Pe} (a student of
Mitchell's) and the already mentioned Lemma 2.4 of \cite{MMR}.

The present note is centered around a series of remarks leading to the conclusion that the convenient and clever trick
presented in \cite{Pe} can be in fact generalized from $R$ to a whole class of countably infinite ultrahomogeneous
structures (that is, \Fr limits), thus yielding the Bergman property for their endomorphism monoids. This conclusion is
reached in the main result of this paper, Theorem \ref{main2}, with the purpose of supplementing the results of
\cite{MMR}. In the following preliminary section, we are going to briefly review the aforementioned trick from
\cite{Pe} and other ingredients needed for our arguments in Section 4. Along the way, in Section 3 we will record an
exact description of \Fr classes whose limits are homomorphism-homogeneous \cite{CN}, accompanied with a number of
examples.

\section{Preliminaries}

\subsection{Strong distortion, coproducts, homomorphism extensions}

The results we are about to revisit are based on another motif connected to the Bergman property, and it traces back to
an old, classical result of W.~Sierpi\'nski, who proved in \cite{S35} that if $X$ is an infinite set, then any
countable set $\{f_i:\ i\geq 0\}$ of self-maps (transformations) $X\to X$ is contained in a 2-generated subsemigroup of
$\T_X$, the semigroup of all self-maps of $X$. Following this landmark example, we say that a semigroup $S$ has
\emph{Sierpi\'nski rank} $n<\omega$ if $n$ is the least positive integer with the property that for any countable
$A\subseteq S$ there exists $s_1,\dots,s_n\in S$ such that $A\subseteq\langle s_1,\dots,s_n\rangle$. If no such $n$
exists, the Sierpi\'nski rank of $S$ is said to be \emph{infinite}. Of course, the Sierpi\'nski rank of a countable
semigroup is simply its \emph{rank}, the minimum size of its generating set, so that the notion is particularly
interesting for uncountable semigroups. So, the result of Sierpi\'nski asserts that $\T_X$ has Sierpi\'nski rank 2.
Some recent results concerning the Sierpi\'nski rank of certain classical transformation semigroups can be found in
\cite{MP}.

A slight modification of this notion produces a convenient method for proving the Bergman property for semigroups.
Namely, in several proofs of the finiteness of the Sierpi\'nski rank for various semigroups it turns out---after
selecting a countable set $A=\{a_i:\ i<\omega\}$ and $s_1,\dots,s_n$ such that $A\subseteq\langle
s_1,\dots,s_n\rangle$---that in the representation $$a_i=\w_i(s_1,\dots,s_n)$$ the length of the word $\w_i$ does not
depend on the particular choice of $a_i$, but that it is determined only by the index $i$. In other words, the
Sierpi\'nski property occurs in some sense in a ``uniform'' way. More formally, call a semigroup $S$ \emph{strongly
distorted} if there exists a sequence of natural numbers $(\ell_n)_{n<\omega}$ and $M<\omega$ such that for any
sequence $(a_n)_{n<\omega}$ of elements of $S$ there exist $s_1,\dots,s_M\in S$ and a sequence of words
$(\w_n)_{n<\omega}$ (over an $M$-letter alphabet) such that $|\w_n|\leq\ell_n$ and $a_n=\w_n(s_1,\dots,s_M)$ for all
$n<\omega$. Here is the result that puts strongly distorted semigroups into the context of the initial motivation of
this paper.

\begin{lem}[{\cite[Lemma 2.4]{MMR}}]\label{sd}
If $S$ is a non-finitely generated and strongly distorted semigroup, then $S$ has the Bergman property.
\end{lem}

Therefore, any strongly distorted uncountable semigroup has the Bergman property. This observation is the link showing
that Lemma 3.10.3 and the proof of Theorem 3.10.4 in \cite{Pe} in fact establish the Bergman property for $\End(R)$.
However, the good thing about the latter theorem is that it is not really about the random graph, as it very easily
admits a generalization that we present here. But first recall the classical category-theoretical notion of a
\emph{coproduct}. If $\{A_i:\ i\in I\}$ is a family of first-order structures belonging to a concrete category
$\mathbf{C}$ (where objects are structures and morphisms are their homomorphisms), then their coproduct (or \emph{free
sum}), denoted by $\coprod\st_{i\in I}A_i$, is a structure $S\in\mathbf{C}$ with the following properties:
\begin{itemize}
\item[(a)] there are embeddings $\iota_i:A_i\to S$ for any $i\in I$;
\item[(b)] for any $B\in\mathbf{C}$ and any homomorphisms $\varphi_i:A_i\to B$, $i\in I$, there exists a unique homomorphism
$\varphi:S\to B$ such that $\varphi\iota_i=\varphi_i$ holds for all $i\in I$.
\end{itemize}
(In this paper, mappings are composed right to left, so that $fg$ is a function for which $fg(x)$ means $f(g(x))$. For
a set $X$, $\mathbf{1}_X$ will always denote the identity mapping on $X$.) Whenever it exists, the coproduct is unique
up to an isomorphism, and it is generated by $\bigcup_{i\in I}\iota_i(A_i)$.

So, here is the ``abstract'' version of Theorem 3.10.4 from \cite{Pe}.

\begin{thm}[Mitchell \cite{JDM}, P\'eresse \cite{Pe}] \label{main0}
Let $A$ be an infinite structure with a sub\-struc\-tu\-re $B$ satisfying the following conditions:
\begin{itemize}
\item[(i)] $B\cong\coprod_{n<\omega}A_n$, where $A_n\cong A$ for each $n<\omega$;
\item[(ii)] any homomorphism $\varphi:B\to A$ can be extended to an endomorphism $\widehat\varphi$ of $A$.
\end{itemize}
Then $\End(A)$, the endomorphism monoid of $A$, is strongly distorted. In addition, the Sierpi\'nski rank of $\End(A)$
is at most 3.
\end{thm}

\beginproof
Let $f_0,f_1,\dots$ be any countable sequence of endomorphisms of $A$. We construct $g_1,g_2,g_3\in\End(A)$ such that
$f_k\in\langle g_1,g_2,g_3\rangle$ for any $k\geq 0$. Also, we will freely assume that each $A_n$ is actually contained
in $B$.

First of all, let $g_1:A\to A_0$ be any isomorphism. Furthermore, let $h_n:A_n\to A_{n+1}$, $n<\omega$, be
isomorphisms. By (i) and the definition of the coproduct, since each $h_n$ maps into $B$, there is a homomorphism
$h:B\to B$ (that is, $h\in\End(B)$) such that $h|_{A_n}=h_n$ for any $n<\omega$. By condition (ii), $h$ can be extended
to an endomorphism of $A$, which we denote by $g_2$. Then for any $n<\omega$ we have $g_2|_{A_n}=h_n$ and
$t_n=g_2^ng_1:A\to A_n$ is an isomorphism.

Now we define the key endomorphism $g_3$, which can be informally thought of as a ``compressed form'' of the sequence
$\{f_k\}_{k\geq 0}$, where each $f_k$ is ``packed up'' into the copy $A_k$ of $A$. Since $t_k^{-1}:A_k\to A$ is an
isomorphism, it follows that $\psi_k=f_kt_k^{-1}$ is a homomorphism of $A_k$ into $A$. Similarly as above, by (i) there
exists a homomorphism $\psi:B\to A$ such that $\psi|_{A_k}=\psi_k$ for each $k\geq 0$. However, by (ii) there is an
extension of $\psi$ to $g_3\in\End(A)$. Note that $g_3|_{A_k}=\psi_k$ holds as well.

It remains to recover $f_k$ from $g_3$. Indeed, since $t_k$ maps into $A_k$,
$$g_3g_2^kg_1=g_3t_k=\psi_kt_k=f_kt_k^{-1}t_k=f_k,$$
as wanted. So, not only $f_k\in\langle g_1,g_2,g_3\rangle$, but we uniformly have that the length of the product
representing $f_k$ is $\ell_k=k+2$.
\endproof

Of course, the coproduct in the category of (simple) graphs is just the disjoint union of the given family of graphs.
Hence, $\coprod_{n<\omega}R$ exists, and, since it is a countably infinite graph, it embeds into $R$. In addition, as
established in Lemma 3.10.3 of \cite{Pe}, $R$ has the remarkable property that for any countable graph $G$ there is an
induced subgraph $G'$ of $R$, isomorphic to $G$, such that any homomorphism $\varphi:G'\to R$ extends to an
endomorphism $\wh\varphi$ of $R$. As $\End(R)$ is known to be uncountable, it follows that it has the Bergman property.

Our main goal here is to see to which extent we can utilize the above theorem in order to cover the cases of some of
the most important infinite structures arising as \Fr limits. The really intriguing condition here is Theorem
\ref{sd}(ii), the possibility of extending a homomorphism from a certain substructure of $\Fl(\C)$. More precisely, we
will be interested in the conditions under which for any $A\in\ol\C$ there exists a substructure $A'$ of $F=\Fl(\C)$
such that $A'\cong A$ and any homomorphism $\varphi:A'\to F$ can be extended to an endomorphism of $F$. These
conditions must be nontrivial, as the next example shows.

\begin{exa}\rm
Let $H_n$ denote the \emph{Henson graph} \cite{He,LW}, that is, the \Fr limit of the class of all finite simple graphs
omitting $K_n$, the complete graph (clique) on $n$ vertices, $n\geq 3$. Now, $H_n$ clearly contains copies of $K_{n-1}$
and $\ol{K_{n-1}}$ and any bijection of vertices $f:\ol{K_{n-1}}\to K_{n-1}$ is a graph homomorphism. However, $f$
cannot be extended to an endomorphism $\wh{f}$ of $H_n$ because there is a vertex $u$ adjacent to each vertex of the
anti-clique $\ol{K_{n-1}}$, so that $\wh{f}(u)$ would be adjacent to each vertex of the clique $K_{n-1}$, which is
impossible by the definition of $H_n$. Therefore, $H_n$ is not homomorphism-homogeneous. In addition, as shown in
\cite{Mu}, every endomorphism of $H_n$ is injective, so not every endomorphism of a (finite) subgraph of $H_n$ can be
extended to a member of $\End(H_n)$ (just take two non-adjacent vertices $u,v$ and map them both into $u$).
\end{exa}

\subsection{The amalgamation property and amalgamated free sums}

An \emph{amalgam} is a quintuple $(A,B,C,f_1,f_2)$ consisting of structures $A,B,C$ together with two embeddings
$f_1:A\to B$ and $f_2:A\to C$. If $A,B,C\in\C$ for some class $\C$, then we have an amalgam \emph{in} $\C$. The
\emph{amalgamation property} for $\C$, mentioned earlier, asserts that any amalgam in $\C$ can be embedded into a
structure $D\in\C$, i.e.\ that there are embeddings $g_1:B\to D$ and $g_2:C\to D$ such that $g_1f_1=g_2f_2$. If $\C$ is
a class of finitely generated structures with the AP (for example, a \Fr class), then it is known that the statement of
the AP extends to non-finitely generated members of $\ol\C$ in the following sense. The proof is rather
straightforward, thus it is omitted.

\begin{lem}\label{AP}
Let $\C$ be a class of finitely generated structures satisfying the amalgamation property. If $(A,B,C,f_1,f_2)$ is an
amalgam such that $A\in\C$ and $B,C\in\ol\C$, then it can be embedded into some structure $D\in\ol\C$.
\end{lem}

In the course of dealing with a \Fr class $\C$ it would be very useful to fix a \emph{canonical} way for embedding an
amalgam $(A,B,C,f_1,f_2)$ such that $A\in\C$ and $B,C\in\ol\C$ into a structure from $\ol\C$. Such possibility is
provided by the standard categorical notion of the \emph{pushout} (see \cite{Mac} for a background in basic category
theory). Recall that if $f:X\to Y$ and $g:X\to Z$ are two morphisms, then their pushout is an object $P$ together with
two morphisms $i_1:Y\to P$ and $i_2:Z\to P$ such that the following diagram commutes:
$$
\xymatrix{%
Y \ar[r]^{i_1} & P\\
X \ar[u]^f \ar[r]_g & Z \ar[u]_{i_2} }
$$
while for any object $Q$ and morphisms $j_1:Y\to Q$ and $j_2:Z\to Q$ such that $j_1f=j_2g$ there exists a unique
morphism $h:P\to Q$ such that $j_1=hi_1$, $j_2=hi_2$:

$$
\xymatrix{&& Q\\
Y \ar[r]^{i_1} \ar@/^2ex/[urr]^{j_1} & P \ar[ru]^h & \\
X \ar[u]^f \ar[r]_g & Z \ar[u]_{i_2} \ar@/_2ex/[ruu]_{j_2} & }
$$
In concrete categories of structures we often consider the case when $f,g$ are embeddings, whence in the presence of
the AP the homomorphisms $i_1,i_2$ must be injective as well. Hence, $P$ can be thought of as the ``smallest''
structure embedding the amalgam $(X,Y,Z,f,g)$.

Accordingly, a structure $P$ will be called the \emph{amalgamated free sum} of $Y$ and $Z$ with respect to $X$ if there
exist embeddings $i_1:Y\to P$ and $i_2:Z\to P$ such that $P$ (with $i_1$ and $i_2$) is the pushout of the amalgam
$(X,Y,Z,f,g)$. If so, we write $P=Y\ast_X Z$. It is easy to check that the amalgamated free sum, if it exists, is
unique up to an isomorphism, and that it is generated by $i_1(Y)\cup i_2(Z)$. We will consider \Fr classes $\C$
satisfying the following condition, which is a rather strong form of the AP:
\begin{itemize}
\item[\cd] For any amalgam $(A,B,C,f_1,f_2)$ such that $A\in\C$, $B,C\in\ol\C$, the amalgamated free sum $B\ast_A C$
exists and belongs to $\ol\C$.
\end{itemize}

For example, it is straightforward to see that the amalgamated free sum of two simple graphs $G_1=(V_1,E_1)$ and
$G_2=(V_2,E_2)$ sharing a common induced subgraph on $V=V_1\cap V_2$ is simply the graph on $V_1\cup V_2$ whose edges
are $E_1\cup E_2$ (in a somewhat simplified form, one may say that the free sum of the amalgam is the amalgam itself).
Similarly, the amalgamated free sum of two posets $(A_1,\leq_1)$ and $(A_2,\leq_2)$ over a common subposet on
$A=A_1\cap A_2$ is the poset $(A_1\cup A_2,\leq)$, where $\leq$ is the transitive closure of the reflexive and
antisymmetric relation $\leq_1\cup\leq _2$ on $A_1\cup A_2$. The existence of amalgamated free sums of semilattices,
distributive lattices and Boolean algebras can be traced from \cite{Gr}; also, the construction of amalgamated free
sums for Abelian groups and vector spaces is a part of algebraic folklore.

\subsection{From amalgamated sums to the construction of \Fr limits}\label{constr}

Let $\C$ be a \Fr class satisfying the condition \cd. Equipped with the construction of the amalgamated free sum, we
first describe a particular extension $A^\star$ for an arbitrary structure $A\in\ol\C$. This is in fact a
generalization of one of the standard constructions of the random graph described in \cite{C1} and an adaptation of the
general approach from \cite{H}.

First of all, recall that a structure $C$ is a \emph{one-point extension} of its substructure $B$ if there is an
element $x\in C\setminus B$ such that $C$ is generated by $B\cup\{x\}$. Trivially, if $B$ is finitely generated, so is
$C$.

Now let $\{(B_i,C_i):\ i<\omega\}$ be the enumeration of all pairs of structures such that $B_i\in\C$ is a finitely
generated substructure of $A$, while $C_i$ is a one-point extension of $B_i$ belonging to $\C$; for each isomorphism
type we take one such extension. We construct a chain of structures $A_i$, $i\geq 0$, as follows. Let $A_0=A$ and
assume that $A_n$ has already been constructed for some $n\geq 0$ such that $A\subseteq A_n\in\ol\C$. Then $B_n$ is a
substructure of $A$ and so of $A_n$, whence $(B_n,A_n,C_n,\mathbf{1}_{B_n},\mathbf{1}_{B_n})$ is an amalgam such that
$B_n,C_n\in\C$ and $A_n\in\ol\C$. The condition \cd\ allows us to define
$$A_{n+1}=A_n\ast_{B_n}C_n.$$
Clearly, $A_{n+1}$ embeds the considered amalgam; therefore, there is no loss of generality in assuming that
$A_n\subseteq A_{n+1}$, so that $A$ is a substructure of $A_{n+1}$. Finally, we let
$$A^\star=\bigcup_{n<\omega}A_n.$$
This construction can be iterated by setting $A^{(0)}=A$ and $A^{(n+1)}=(A^{(n)})^\star$ for all $n\geq 0$. Let
$$\F(A)=\bigcup_{n<\omega}A^{(n)},$$
which is an extension of $A$. Clearly, any finitely generated substructure of $\F(A)$ must belong to some $A^{(m)}$,
and since $A^{(m)}\in\ol\C$ the finitely generated structure in question belongs to $\C$; hence, $\F(A)\in\ol\C$.

\begin{pro}\label{flim}
Let $\C$ be a \Fr class satisfying \cd. For any $A\in\ol\C$, the structure $\F(A)$ constructed as above is isomorphic
to the \Fr limit of $\C$.
\end{pro}

\beginproof
This is a consequence of another well-known result of \Fr \cite{F2} (see also \cite[Lemma 6.1.3]{H}): namely, it
suffices to prove that $\F(A)$ \emph{realizes all one-point extensions in} $\C$ (effectively, this says that $\Fl(\C)$
is the unique existentially closed structure in $\ol\C$). This means that for each finitely generated substructure $B$
of $\F(A)$ and its one-point extension $C\in\C$ there should be an embedding $f:C\to\F(A)$ such that
$f|_B=\mathbf{1}_B$.

However, this is very easy to check. Namely, as already remarked, $B$ must be a substructure of $A^{(m)}$ for some
$m\geq 0$. Hence, if $\{(B_i^{(m)},C_i^{(m)}):\ i<\omega\}$ is the enumeration of pairs of structures required for the
construction of $A^{(m+1)}=(A^{(m)})^\star$, then $B=B_j^{(m)}$ and $C\cong C_j^{(m)}$ for some $j$, with an
isomorphism $f_0:C\to C_j^{(m)}$ such that $f_0|_B=\mathbf{1}_B$. So, in the course of producing $A_{j+1}^{(m)}$ from
$A_j^{(m)}$ we embed the amalgam $(B,A_j^{(m)},C_j^{(m)},\mathbf{1}_B,\mathbf{1}_B)$ into $A_{j+1}^{(m)}$, and so into
$A^{(m+1)}$. Hence, there exists an embedding $f_1:C_j^{(m)}\to A_{j+1}^{(m)}\subseteq\F(A)$ such that
$\mathbf{1}_Bf_1=\mathbf{1}_B\mathbf{1}_{A_j^{(m)}}$ and so $f_1|_B=\mathbf{1}_B$. Now $f=f_1f_0:C\to\F(A)$ is the
required embedding, since $f|_B=\mathbf{1}_B$.
\endproof

For example, for an arbitrary countable graph $G=(V,E)$, the graph $G^\star$ is obtained by adjoining a vertex $u_A$
for each finite subset $A\subseteq V$ such that $u_A$ is joined by an edge to $v\in V$ if and only if $v\in A$. By
iterating this construction, we build a countable graph $R(G)$ ``around'' the initial graph $G$. As remarked in
\cite{C1}, regardless of the choice of $G$, we always end up with $R(G)\cong R$, the countably infinite random graph.

\begin{rem}\rm
Of course, Proposition \ref{flim} holds in greater generality, for arbitrary \Fr classes. The construction of the
structure $\F(A)$ should be amended in the more general case so that $A_{n+1}$ is selected to be an arbitrary structure
from $\ol\C$ embedding the amalgam $(B_n,A_n,C_n,\mathbf{1}_{B_n},\mathbf{1}_{B_n})$ and containing $A_n$ as a
substructure. (Such $A_{n+1}$ exists by Lemma \ref{AP}.) However, we will need the condition \cd\ and the previous more
specific construction of $\F(A)$ in Section 4, where we present our main arguments; in particular, this will be crucial
in the proof of Theorem \ref{main1}.
\end{rem}

\section{Homomorphism-homogeneous \Fr limits}

In an attempt to put the notion of ultrahomogeneity into a more general setting of arbitrary homomorphisms of
first-order structures, Cameron and Ne\v set\v ril introduced in \cite{CN} the property of
\emph{homomorphism-homogeneity}. Namely, a structure $A$ is homomorphism-homogeneous if any homomorphism $B\to C$
between its finitely generated substructures can be extended to an endomorphism of $A$. Recent results concerning
characterizations of this property in various classes of structures include \cite{CL,DM,IMR,M,RSch}. In this section,
we make a brief pause towards our aim to record a condition equivalent to homomorphism-homogeneity of a \Fr limit.

To this end, we introduce yet another property that a \Fr class $\C$ may or may not satisfy. We say that $\C$ satisfies
the \emph{one-point homomorphism extension property (1PHEP)} if for any $B,B',C\in\C$ such that $C$ is a one-point
extension of $B$, $C=\langle B\cup\{x\}\rangle$, any surjective homomorphism $\varphi:B\to B'$ can be extended to a
homomorphism $\varphi':C\to C'$ for some $C'\in\C$ containing $B'$. If we require that $\varphi'$ is surjective as
well, then it is clear that either $C'=B'$, or $C'$ is a one-point extension of $B'$, namely $C'=\langle
B'\cup\{\varphi'(x)\}\rangle$.

\begin{rem}\rm
It is quite easy to show that for any class of finitely generated structures of a fixed signature, the 1PHEP is
equivalent to the seemingly more general \emph{homo-amalgamation property (HAP)}, which, even though it is not
explicitly formulated, transpires from the treatment in Section 4 of \cite{CN}. Namely, the HAP is the assertion that
for any $A,B_1,B_2\in\C$, any homomorphism $\varphi:A\to B_1$ and any embedding $f:A\to B_2$ there is a structure
$D\in\C$, an embedding $f':B_1\to D$ and a homomorphism $\varphi':B_2\to D$ such that $f'\varphi=\varphi'f$. However,
the more specific form of the 1PHEP might be slightly easier to check, as the following examples show.
\end{rem}

\begin{exa}\rm
The class of all finite simple graphs has the 1PHEP. Indeed, let $\varphi:G\to H$ be a surjective graph homomorphism,
and let $G'$ be a graph obtained from $G$ by adjoining a new vertex $x$ (and some new edges involving $x$). Construct a
new graph $H'$ obtained by adjoining a new vertex $x'$ to $H$, while for $v\in V(H)$ we set that $(x',v)\in E(H')$ if
and only if $(x,u)\in E(G')$ for some $u\in V(G)$ such that $\varphi(u)=v$. Then it is easily verified that
$\varphi':G'\to H'$ obtained by extending $\varphi$ by $\varphi'(x)=x'$ is a (surjective) graph homomorphism.
\end{exa}

\begin{exa}\rm
We have already seen that the \Fr class of all $K_n$-free finite simple graphs fails to satisfy the 1PHEP: a bijection
from the vertices of an anti-clique of size $n-1$ to a clique of the same size cannot be extended within the considered
class to a vertex adjacent to all vertices of the anti-clique.

On the other hand, the ``complementary'' \Fr class to the above one, that of all $\ol{K_n}$-free finite simple graphs
has the 1PHEP: it is quite straightforward to check that the construction from the previous example will work for this
class as well.
\end{exa}

\begin{exa}\rm
The class of all finite posets has the 1PHEP. To see this, let $B$ be a finite poset, $C=B\cup\{x\}$ its one-point
extension, and $\varphi:B\to B'$ an order-preserving map (a poset homomorphism). Let $$L=\{b\in B:\
b<x\}\quad\mbox{and}\quad U=\{b\in B:\ x<b\}.$$ Since $\ell<u$ holds for any $\ell\in L$ and $u\in U$ we have
$\varphi(\ell)\leq\varphi(u)$. So, $\varphi(L)\cap\varphi(U)$ is either empty, or a singleton. In the former case,
define an exten\-sion $C'$ of $B'$ by ``inserting'' a new element $y$ between $L'=\varphi(L)$ and $U'=\varphi(U)$; this
is possible as $\ell'<u'$ holds for any $\ell'\in L'$, $u'\in U'$. It is now a routine to check that the mapping
$\varphi'$ such that $\varphi'|_B=\varphi$ and $\varphi'(x)=y$ is a poset homomorphism $C\to C'$. If, however,
$\varphi(L)\cap\varphi(U)=\{x'\}$ then extend $\varphi$ to $\varphi':C\to B'$ by defining $\varphi'(x)=x'$; once again,
$\varphi'$ turns out to be a homomorphism.
\end{exa}

Recall that \emph{metric spaces} can be viewed as first-order structures over an uncountable language consisting of
binary relational symbols indexed by the non-negative reals such that $(x,y)\in R_\alpha$ ($\alpha\in\mathbb{R}_0^+$)
if and only if $d(x,y)\leq\alpha$. (Of course, we may as well restrict ourselves to metric spaces with rational
distances, thus obtaining a countable signature for such structures.) From such a point of view, homomorphisms of
metric spaces are just \emph{non-expanding functions} $\varphi$ so that we have
$$d(\varphi(x),\varphi(y))\leq d(x,y)$$
for any $x,y$. Then, naturally, the notion of an automorphism coincides with that of an \emph{isometry}, a
distance-preserving permutation.

\begin{lem}\label{metr}
The class of finite metric spaces has the 1PHEP. The same applies to (the \Fr class of) finite metric spaces with
rational distances.
\end{lem}

\beginproof
Let $\varphi:M\to M'$ be a surjective homomorphism of finite metric spaces, and let $M_1$ be a one-point extension of
$M$, with $y$ being the new point. Our aim is to prove that there exists a metric space $M_1'=M'\cup\{y'\}$, a
one-point extension of $M'$, such that for each $x\in M$ we have $d(y',\varphi(x))\leq d(y,x)$. Then we can extend
$\varphi$ to a homomorphism $\wh\varphi:M_1\to M_1'$ by defining $\wh\varphi(y)=y'$.

Let $M=\{x_i:\ i<n\}$.

First of all, we are going to consider the special case when $\varphi$, the initial homomorphism, is a bijection. We
are looking for a sequence of positive real numbers $\dd_i$, $i<n$, such that a (hypothetical) point $y'$ with
$d(y',\varphi(x_i))=\dd_i$ for all $i<n$ satisfies all the triangle inequalities with the already existing points of
$M'$. In other words, the required conditions are:
\begin{itemize}
\item[(1)] $\dd_i+\dd_j \geq d(\varphi(x_i),\varphi(x_j))$,
\item[(2)] $\dd_i+d(\varphi(x_i),\varphi(x_j)) \geq \dd_j$,
\end{itemize}
with $i,j<n$, $i\neq j$, in both cases. In addition, we need the third condition
\begin{itemize}
\item[(3)] $\dd_i\leq d(y,x_i)$ for all $i<n$.
\end{itemize}

There is no loss of generality in assuming that $d(y,x_0)\leq d(y,x_1)\leq\dots\leq d(y,x_{n-1})$. Now consider the
sequence defined by $\dd_0=d(y,x_0)$ and
$$\dd_i=\min_{0\leq k<i}\{d(y,x_i), d(y,x_k)+d(\varphi(x_k),\varphi(x_i))\}$$
for $0<i<n$. We claim that these numbers constitute a solution of the system of inequalities (1)--(3) above. Indeed,
the condition (3) is immediately satisfied. For (1), we distinguish three subcases. If
$\dd_i=d(y,x_k)+d(\varphi(x_k),\varphi(x_i))$ and $\dd_j=d(y,x_m)+d(\varphi(x_m),\varphi(x_j))$ for some $k<i$ and
$m<j$, then
\begin{align*}
\dd_i+\dd_j&=d(\varphi(x_i),\varphi(x_k))+d(y,x_k)+d(y,x_m)+d(\varphi(x_m),\varphi(x_j))\\
&\geq d(\varphi(x_i),\varphi(x_k))+d(x_k,x_m)+d(\varphi(x_m),\varphi(x_j))\\
&\geq d(\varphi(x_i),\varphi(x_k))+d(\varphi(x_k),\varphi(x_m))+d(\varphi(x_m),\varphi(x_j))\\
&\geq d(\varphi(x_i),\varphi(x_j)),
\end{align*}
since $M,M'$ are metric spaces and $\varphi$ is non-expanding. On the other hand, if $\dd_i=d(y,x_i)$ and
$\dd_j=d(y,x_j)$, then
$$\dd_i+\dd_j=d(y,x_i)+d(y,x_j)\geq d(x_i,x_j)\geq d(\varphi(x_i),\varphi(x_j)).$$
Finally, if $\dd_i=d(y,x_k)+d(\varphi(x_k),\varphi(x_i))$ for some $k<i$ and $\dd_j=d(y,x_j)$ (the symmetric case is
analogous), then
\begin{align*}
\dd_i+\dd_j&=d(y,x_j)+d(y,x_k)+d(\varphi(x_k),\varphi(x_i))\\
&\geq d(x_j,x_k)+d(\varphi(x_k),\varphi(x_i))\\
&\geq d(\varphi(x_j),\varphi(x_k))+d(\varphi(x_k),\varphi(x_i))\\
&\geq d(\varphi(x_j),\varphi(x_i))=d(\varphi(x_i),\varphi(x_j)).
\end{align*}
Concerning (2), assume first that $i<j$. Then if $\dd_i=d(y,x_i)$ we have
$$\dd_i+d(\varphi(x_i),\varphi(x_j))=d(y,x_i)+d(\varphi(x_i),\varphi(x_j))\geq \dd_j$$ by the definition of $\dd_j$; if,
however, $\dd_i=d(y,x_k)+d(\varphi(x_k),\varphi(x_i))$ for some $k<i$ then
\begin{align*}
\dd_i+d(\varphi(x_i),\varphi(x_j))&=d(y,x_k)+d(\varphi(x_k),\varphi(x_i))+d(\varphi(x_i),\varphi(x_j))\\
&\geq d(y,x_k)+d(\varphi(x_k),\varphi(x_j))\geq \dd_j,
\end{align*}
as $k<i<j$. So, it remains to discuss the possibility $i>j$. If it happens that
$\dd_i=d(y,x_k)+d(\varphi(x_k),\varphi(x_i))$ for some $k<i$, and, in addition, we have $k<j$ as well, then
$\dd_i+d(\varphi(x_i),\varphi(x_j))\geq \dd_j$ holds by the identical argument as in the previous displayed chain of
equalities and inequations. Otherwise, either $k\geq j$, or $\dd_i=d(y,x_i)$, both cases implying $\dd_i\geq
d(y,x_k)\geq d(y,x_j)\geq \dd_j$, so (2) holds. This completes the case when $\varphi$ is injective, since $M'_1$ is
obtained by adjoining a point $y'$ to $M'$ such that $d(y',\varphi(x_i))=\dd_i$ for all $i<n$.

Turning to the general case, when $\varphi$ is not necessarily a bijection, for any $z\in M'$ choose a point
$\mathbf{x}_z\in\varphi^{-1}(z)\subseteq M$ whose distance to $y$ is minimal among all elements of $\varphi^{-1}(z)$
(that is, we have $d(y,\mathbf{x}_z)\leq d(y,x)$ for all $x\in M$ such that $\varphi(x)=z$). Let $M_0=\{\mathbf{x}_z:\
z\in M'\}$. Now $\varphi|_{M_0}:M_0\to M'$ is a bijective homomorphism of finite metric spaces, so by the previous
considerations it follows that there is a one-point extension $M'_1$ of $M'$ and a homomorphism $\psi:M_0\cup\{y\}\to
M'_1$ extending $\varphi|_{M_0}$. But then $\wh\varphi=\psi\cup\varphi$ is the required extension of $\varphi$, since
for any $x\in M$ we have
$$d(\wh\varphi(y),\wh\varphi(x))=d(y',\varphi(x))=d(\psi(y),\psi(\mathbf{x}_{\varphi(x)}))\leq
d(y,\mathbf{x}_{\varphi(x)})\leq d(y,x),$$ as wanted. It remains to note that if all
$d(y,x_i),d(x_i,x_j),d(\varphi(x_i),\varphi(x_j))$ are rational numbers, so are all $\dd_i$, thus the second part of
the assertion follows, too.
\endproof

The 1PHEP occurs in algebraic structures as well, where it is intimately related to the \emph{congruence extension
property (CEP)}, see \cite{Gr}. Recall that an algebra $A$ has the CEP if for any subalgebra $B$ of $A$ and any
congruence $\rho$ of $B$ there exists a congruence $\theta$ of $A$ whose restriction to $B$ is precisely $\rho$, that
is, $\theta\cap(B\times B)=\rho$.

\begin{lem}
Let $\C$ be a \Fr class of algebras closed for taking homomorphic images. If all members of $\C$ have the CEP, then
$\C$ has the 1PHEP.
\end{lem}

\beginproof
Let $C=\langle B\cup\{x\}\rangle$ and let $\varphi:B\to B'$ be a surjective homomorphism, where $B,B',C\in\C$. Then
$\rho=\krn\varphi$ is a congruence of $B$ (such that $B/\mathrm{ker}\,\varphi \cong B'$), so by the CEP there exists a
congruence $\theta$ of $C$ such that $\theta\cap(B\times B)=\krn\varphi$. Now consider the natural homomorphism
$\nu_\theta:C\to C/\theta$. The image of $B$, $\nu_\theta(B)$, is isomorphic to $B/(\theta\cap(B\times B))$, which is
by the given conditions isomorphic to $B'$. Therefore, $C/\theta\in\C$ can be considered as an extension of $B'$,
whence $\nu_\theta$ is an extension of $\varphi$.
\endproof

By invoking the fact that the varieties of semilattices, distributive lattices, Boolean algebras and vector spaces over
a given field $\mathbb{F}$ all possess the CEP, we obtain the following conclusion.

\begin{cor}\label{1phep-ex}
Each of the \Fr classes of all finite semilattices, all finite distributive lattices, all finite Boolean algebras and
all finite-dimensional vector spaces over a field $\mathbb{F}$ have the 1PHEP.
\end{cor}

Now we provide a characterization of homomorphism-homogeneous \Fr limits. It reduces a property of the intricate
structure of such a limit to a ``local'' property of finitely generated structures from $\C$ which usually have much
more transparent features. A related result is contained in \cite[Proposition 4.1]{CN}.

\begin{pro}\label{hom-hom}
Let $\C$ be a \Fr class. Then the \Fr limit of $\C$ is homomor\-phism-homogeneous if and only if $\C$ has the 1PHEP.
\end{pro}

\beginproof
Throughout the proof, let $F=\Fl(\C)$.

($\Ra$) Let $B,B',C\in\C$ be such that $C$ is a one-point extension (or any extension of finite relative rank, for that
matter) of $B$, and let $\xi:B\to B'$ be a surjective homomorphism. By the properties of the \Fr limit, there is no
loss of generality if we assume that $B,B',C$ are in fact substructures of $F$. However, by the given conditions then
there is a $\wh\xi\in\End(F)$ extending $\xi$, so that $\xi'=\wh\xi|_{C}:C\to\wh\xi(C)$ is the homomorphism required by
the 1PHEP.

($\Leftarrow$) Let $A$ be a finitely generated substructure of $F$, while $\varphi:A\to F$ is a homomorphism. Then,
since $F$ is countable, there exists a chain $\{F_i:\ i<\omega\}$ of finitely generated substructures of $F$ such that
$F_0=A$, $F_{i+1}$ is a one-point extension of $F_i$ for each $i\geq 0$, and $F=\bigcup_{i<\omega}F_i$. We construct by
induction a chain of homomorphisms $\varphi_i:F_i\to F$ starting with $\varphi_0=\varphi$. By the 1PHEP, given
$\varphi_j$ for some $j\geq 0$, there exists a finitely generated structure $B'_{j+1}\in\C$, which is an extension of
$B_j=\varphi_j(F_j)$, and a homomorphism $\psi_{j+1}:F_{j+1}\to B'_{j+1}$ that extends $\varphi_j$ (i.e.\
$\psi_{j+1}|_{B_j}=\varphi_j$). Now since $F$ is the \Fr limit of $\C$, there exists an embedding $f_{j+1}:B'_{j+1}\to
F$ which is the identity mapping on $B_j$; define $B_{j+1}=f_{j+1}(B'_{j+1})\subseteq F$. Whence,
$\varphi_{j+1}=f_{j+1}\psi_{j+1}$ is a homomorphism $F_{j+1}\to F$ that extends $\varphi_j$. It remains to define
$$\wh\varphi=\bigcup_{i<\omega}\varphi_i$$
to obtain an endomorphism of $F$ that extends $\varphi$.
\endproof

By combining the previous proposition, Corollary \ref{1phep-ex}, Lemma \ref{metr} and the examples that precede it, we
arrive at the following result.

\begin{cor}
Each of the following \Fr limits is homomorphism-homogeneous: $R$, $\ol{H_n}$ for all $n\geq 3$,
$\mathbb{U}_\mathbb{Q}$, $\mathbb{P}$, $\Omega$, $\mathbb{D}$, $\mathbb{A}$, and $V_\infty^\mathbb{F}$, the
$\aleph_0$-dimensional vector space over a field $\mathbb{F}$.
\end{cor}

\begin{rem}\label{Ury}\rm
Based on Lemma \ref{metr} and an analogous approach as in Proposition \ref{hom-hom}, it is now quite easy to prove that
the \emph{Urysohn space} $\mathbb{U}$ \cite{Ur}, the completion of $\mathbb{U}_\mathbb{Q}$, is homomorphism-homogeneous
as well. Namely, if $X\subseteq\mathbb{U}$ is a finite metric space and $\varphi:X\to\mathbb{U}$ is a homomorphism,
then one can select a countable dense subspace $Y\subseteq\mathbb{U}$ isometric to $\mathbb{U}_\mathbb{Q}$ and use
Lemma \ref{metr} to obtain a homomorphism $\varphi':X\cup Y\to\mathbb{U}$ extending $\varphi$. Now it remains to remark
that: (a) $\mathbb{U}$ is the completion of $X\cup Y$, and (b) every homomorphism of metric spaces is a uniformly
continuous mapping (since it is in fact a Lipschitz function with constant $1$), whence the properties of the
completion of a metric space yield an endomorphism $\wh\varphi$ of $\mathbb{U}$ extending $\varphi$.
\end{rem}

\section{Homomorphism extensions and the Bergman property}

We start immediately with a condition ensuring that an instance of a \Fr limit---as constructed in Subsection
\ref{constr}---satisfies the condition (ii) from Theorem \ref{main0}, related to the possibility of extending a partial
endomorphism of such a limit.

\begin{thm}\label{main1}
Let $\C$ be a \Fr class satisfying \cd\ and the 1PHEP, and let $A\in\ol\C$. Then every homomorphism $\varphi:A\to\F(A)$
can be extended to a $\wh\varphi\in\End(\F(A))$.
\end{thm}

\beginproof
Let $\varphi:A\to\F(A)$ be any homomorphism. Our aim is to obtain a sequence of homomorphisms
$\varphi^{(0)}=\varphi\subseteq\varphi^{(1)}\subseteq\varphi^{(2)}\subseteq\dots$, where $\varphi^{(n)}:A^{(n)}\to
\F(A)$, whence
$$\wh\varphi=\bigcup_{n<\omega}\varphi^{(n)}$$
will be the desired endomorphism of $\F(A)$. Therefore, we start with the assumption that the required sequence has
already been constructed up to $\varphi^{(n)}$ for some $n\geq 0$.

In addition, recall that $A^{(n+1)}$ has been obtained from $A^{(n)}$ by successive amalgamations of all possible (up
to isomorphism) one-point $\C$-extensions
$$\{(B_i^{(n)},C_i^{(n)}):\ i<\omega\}$$ of finitely generated substructures of $A^{(n)}$.
This results in a sequence of structures $A_0^{(n)}=A^{(n)}\subseteq A_1^{(n)}\subseteq\dots$ whose limit is
$(A^{(n)})^\star=A^{(n+1)}$. Accordingly, we construct a tower of homomorphisms $\varphi_i^{(n)}:A_i^{(n)}\to\F(A)$,
$i\geq 0$, as follows, starting with $\varphi_0^{(n)}=\varphi^{(n)}$ and assuming that $\varphi_k^{(n)}$ has already
been constructed.

Now, since $\C$ satisfies \cd, we know that $A_{k+1}^{(n)}$ is obtained as the amalgamated free sum of
$$(B_k^{(n)},A_k^{(n)},C_k^{(n)},\mathbf{1}_{B_k^{(n)}},\mathbf{1}_{B_k^{(n)}}).$$ For
brevity, denote $B'=\varphi^{(n)}(B_k^{(n)})$ and consider the homomorphism between finitely generated $\C$-structures
$\phi=\varphi^{(n)}|_{B_k^{(n)}}:B_k^{(n)}\to B'$. Since $B'$ is finitely generated, there exists an index $p<\omega$
such that $B'\subseteq A^{(p)}$. By the 1PHEP, there exist a structure $C'\in\C$---that is either $B'$, or its
one-point extension---and a surjective homomorphism $\varepsilon:C_k^{(n)}\to C'$ agreeing with $\varphi^{(n)}$ (that
is, with $\phi$) on $B_k^{(n)}$. Moreover, if $C'\neq B'$, then the extension $(B',C')$ can be identified (up to
isomorphism) with $(B_j^{(p+1)},C_j^{(p+1)})$ for some $j$. In any case, we may assume that $\varepsilon(C_k^{(n)})$ is
a (finitely generated) substructure of $A^{(p+1)}$.

What we have right now is depicted in the following diagram:
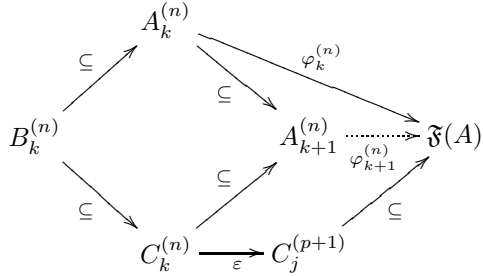
\begin{displaymath}
\xymatrix{%
& A_k^{(n)} \ar[rd]_{\subseteq} \ar[rrd]^{\varphi_k^{(n)}} & & \\
B_k^{(n)} \ar[ru]^{\subseteq} \ar[rd]_{\subseteq} & & A_{k+1}^{(n)}\ar@{.>}[r]_{\varphi_{k+1}^{(n)}\ \ } & \F(A) \\
& C_k^{(n)} \ar[ru]^{\subseteq} \ar[r]_{\varepsilon} & C_j^{(p+1)} \ar[ru]_{\subseteq} & }
\end{displaymath}
By \cd\ and the choice of $A_{k+1}^{(n)}$, there exist a homomorphism $\varphi_{k+1}^{(n)}:A_{k+1}^{(n)}\to \F(A)$ (see
the ``dotted'' arrow) completing the above diagram to a commutative one. In particular, $\varphi_{k+1}^{(n)}$ is an
extension of $\varphi_k^{(n)}$. Finally,
$$\varphi^{(n+1)}=\bigcup_{i<\omega}\varphi_i^{(n)}$$
is a homomorphism $A^{(n+1)}\to \F(A)$, and so we are done.
\endproof

The combination of Lemma \ref{sd}, Theorem \ref{main0} and the previous theorem immediately yields the principal result
of this paper.

\begin{thm}\label{main2}
Let $\C$ be a \Fr class satisfying \cd\ and the 1PHEP, and let $F=\Fl(\C)$. If the coproduct of countably infinitely
many copies of $F$ exists and belongs to $\ol\C$, then $\End(F)$ is strongly distorted and its Sierpi\'nski rank is at
most 3. If, in addition, $\End(F)$ is not finitely generated, then it has the Bergman property.
\end{thm}

Note that one may equivalently replace `not finitely generated' in the above theorem by `uncountable': indeed, if
$\End(F)$ would be countable, then if would be necessarily finitely generated, because of the finite Sierpi\'nski rank.
However, it is again Theorem \ref{main1} that admits to easily establish $|\End(F)|>\aleph_0$, for limits $F$ of
certain \Fr classes $\C$, as the latter inequality will follow from the existence of a structure $A\in\ol\C$ such that
$|\End(A)|>\aleph_0$. We record the following remark, which is of independent interest as well. We say that a semigroup
$T$ \emph{divides} a semigroup $S$ if $T$ is a homomorphic image of a subsemigroup of $S$.

\begin{lem}\label{div}
Let $\C$ be a \Fr class with \cd\ and the 1PHEP, and let $F=\Fl(\C)$. Then for any structure $A\in\ol\C$ we have that
$\End(A)$ divides $\End(F)$.
\end{lem}

\beginproof
By Theorem \ref{main1} we have that $F$ contains an isomorphic copy $A'$ of $A$ such that any endomorphism of $A'$ can
be extended to an endomorphism of $F$. Now consider only those endomorphisms $f$ of $F$ that induce (by restriction) an
endomorphism of $A'$, that is, $f(A')\subseteq A'$. Such endomorphisms form a subsemigroup $S$ of $\End(F)$. Now for
$f,g\in S$ let $(f,g)\in\rho$ if and only if $f|_{A'}=g|_{A'}$. It is easily seen that $\rho$ is a congruence on $S$;
by the given conditions, $S/\rho\cong\End(A')\cong\End(A)$. \hfill $\Box$

\begin{cor}\label{berg}
The endomorphism monoid of any of $R$, $\mathbb{P}$, $\Omega$, $\mathbb{D}$, $\mathbb{A}$ and $V_\infty^\mathbb{F}$ has
the Bergman property, and is divided by $\T_{\aleph_0}$.
\end{cor}

\beginproof
For $R$, the countably infinite anti-clique works as $A$ in Lemma \ref{div}, since now any self-map is an endomorphism
of $A$. Similarly, the countably infinite anti-chain $A$ shows the required assertions about $\End(\mathbb{P})$.
Finally, for limits of \Fr classes of algebras it suffices to note that any self-map on $X$, the set of free generators
of the corresponding free algebra $F(X)$, induces an endomorphism of $F(X)$, whence we let $X$ to be countably
infinite.
\endproof

For various reasons, a number of \Fr classes and their corresponding limits remain outside the scope of this approach.
As we have seen, some of them, such as the finite $K_n$-free simple graphs, fail to have the 1PHEP. Other classes, such
as the finite $\ol{K_n}$-free simple graphs and the finite linear orders have the 1PHEP, and they even have both
coproducts and amalgamated free sums in certain broader concrete categories (of simple graphs and posets,
respectively), but these sums fail to be $\ol{K_n}$-free in the former case, or linearly ordered in the latter.
Finally, some structures simply do not have coproducts and/or amalgamated free sums. For example, there seems to be no
meaningful notion of a coproduct for (rational) metric spaces. This stems from the fact that when we are given a finite
metric space $M$ and we wish to add a new point $x$, then the set of possible vectors of distances of $x$ to the
existing elements of $M$ is in general an unbounded subset of $\mathbb{R}^m$, where $m=|M|$ (or $\mathbb{Q}^m$, if we
go for rational distances), thus rendering impossible the choice of the ``farthest point'' from $M$---something which
would be required should the coproduct of $M$ and a singleton space exist. So, since the linear order $\mathbb{Q}$ and
the universal rational metric space $\mathbb{U}_\mathbb{Q}$ are historically the oldest examples of \Fr limits, it is
natural to ask the following questions.

\begin{prb}\rm
Does the monoid of all order-preserving self-maps of $\mathbb{Q}$ have the Bergman property? More generally, what is
the case with doubly homogeneous linear orders \cite{DH}\,?
\end{prb}

\begin{prb}\rm
Does the endomorphism monoid of $\mathbb{U}_\mathbb{Q}$ have the Bergman property? What about the monoid of all
Lipschitz functions of $\mathbb{U}_\mathbb{Q}$?
\end{prb}

\begin{prb}\rm
Do the endomorphism monoids of ultrahomogeneous graphs $H_n$ and $\ol{H_n}$, $n\geq 3$, have the Bergman property?
\end{prb}

Also, the following tantalizing problem arises.

\begin{prb}\rm
Determine the Sierpi\'nski rank of $\End(R)$ exactly: is it 2 or 3\,? The same question applies to any \Fr limit
mentioned in Corollary \ref{berg}.
\end{prb}

\noindent\textbf{Acknowledgements.} The author is indebted a great deal to an anonymous referee whose thorough reading
of the initial manuscript substantially improved the presentation of the results. I am also very grateful to James D.\
Mitchell, Nik Ru\v skuc (Uni\-ver\-sity of St Andrews) and Dragan Ma\v sulovi\'c (University of Novi Sad) for valuable
discussions and correspondence concerning the topic of this note.

% REFERENCES

\normalsize

\begin{flushleft}
Department of Mathematics and Informatics, University of Novi Sad, Trg Do\-si\-te\-ja Obrado\-vi\-\'ca 4, 21101 Novi
Sad, Serbia\\
\texttt{dockie@dmi.uns.ac.rs}
\end{flushleft}

\end{document}